\documentclass[11pt]{article}
\usepackage{amssymb,psfig,latexsym}

\newtheorem{theorem}{Theorem}[section]
\newtheorem{lemma}{Lemma}[section]

\newcommand{\bsq}{\vrule height .9ex width .8ex depth -.1ex}

\newcommand{\hsp}{\hspace{\parindent}}

\newcommand{\ep}{\epsilon}

\newcommand{\RR}{{\Bbb R}}
\newcommand{\QQ}{{\Bbb Q}}

\newcommand{\ZZ}{{\Bbb Z}}

\newcommand{\sE}{{\cal E}}

\newcommand{\sO}{{\cal O}}

\newcommand{\sR}{{\cal R}}
\newcommand{\sS}{{\cal S}}

\newcommand{\bv}{{\bf v}}

\newcommand{\bx}{{\bf x}}

\newcommand{\beql}[1]{\begin{equation}\label{#1}}
\newcommand{\eqn}[1]{(\ref{#1})}
\newcommand{\eeq}{\end{equation}}

\catcode`\@=11
\renewcommand{\section}{
        \setcounter{equation}{0}
        \@startsection {section}{1}{\z@}{-3.5ex plus -1ex minus
        -.2ex}{2.3ex plus .2ex}{\large\bf}
        }
\catcode`@=12
\makeatletter
\def\eqalignno#1{\displ@y \ta {\bf s} kip\@centering
  \halign to\displaywidth{\hfil$\@lign\displaystyle{##}$\ta {\bf s} kip\z@skip
    & $\@lign\displaystyle{{}##}$\hfil\ta {\bf s} kip\@centering
    & \llap{$\@lign##$}\ta {\bf s} kip\z@skip\crcr
    #1\crcr}}
\makeatother

\makeatletter
\def\@sect#1#2#3#4#5#6[#7]#8{\ifnum #2>\c@secnumdepth
     \def\@svsec{}\else 
     \refstepcounter{#1}\edef\@svsec{\csname the#1\endcsname.\hskip .75em }\fi
     \@tempskipa #5\relax
      \ifdim \@tempskipa>\z@ 
        \begingroup #6\relax
          \@hangfrom{\hskip #3\relax\@svsec}{\interlinepenalty \@M #8\par}%
        \endgroup
       \csname #1mark\endcsname{#7}\addcontentsline
         {toc}{#1}{\ifnum #2>\c@secnumdepth \else
                      \protect\numberline{\csname the#1\endcsname}\fi
                    #7}\else
        \def\@svsechd{#6\hskip #3\@svsec #8\csname #1mark\endcsname
                      {#7}\addcontentsline
                           {toc}{#1}{\ifnum #2>\c@secnumdepth \else
                             \protect\numberline{\csname the#1\endcsname}\fi
                       #7}}\fi
     \@xsect{#5}}
\def\@theorem#1#2{\it \trivlist \item[\hskip \labelsep{\bf #1\ #2.}]}
\makeatother

\begin{document}
\begin{center}
{\large {\bf Dynamics of a Family of Piecewise-Linear 
Area-Preserving Plane Maps}} \\

{\large {\bf III. Cantor Set Spectra}}\\ \bigskip

{\large {\em Jeffrey C. Lagarias}} \\ \smallskip
Department of Mathematics  \\
University of Michigan\\
Ann Arbor, MI 48109-1043 \\
{\tt email:}  {\tt lagarias@umich.edu}\bigskip \\

{\large {\em Eric Rains}} \\ \smallskip
Dept. of Mathematics \\
University of California, Davis\\
Davis, CA 95616-8633  \\
{\tt email:}  {\tt rains@math.davis.edu}\bigskip \\
\vspace*{2\baselineskip}
(August 3, 2006 version) \\

\vspace*{1.5\baselineskip}
{\bf ABSTRACT}
\end{center}

This paper studies the behavior under iteration of 
the  maps $T_{ab}(x, y) = (F_{ab}(x) - y, x)$ of the
plane $\RR^2$, in which $F_{ab}(x) = ax$ if $ x \ge 0$ and $bx$ if $x < 0$.
These maps are area-preserving homeomorphisms of $\RR^2$ that map
rays from the origin into rays from the origin.
Orbits of the map correspond to solutions of the 
nonlinear difference equation
$x_{n+2}= 1/2(a-b)|x_{n+1}| + 1/2(a+b)x_{n+1} - x_n.$
This difference equation can be rewritten in an eigenvalue form 
for a  nonlinear difference operator of Schr\"{o}dinger type
$-x_{n+2}+ 2x_{n+1}-x_{n} + V_{\mu}(x_{n+1})x_{n+1}= Ex_{n+1},$
in which $\mu = \frac{1}{2}(a - b)$ is fixed, and
$V_{\mu}(x)= \mu (sgn(x))$ is an antisymmetric step function
potential, and the energy $E= 2 - \frac{1}{2}(a + b)$.
We study the set $\Omega_{SB}$ of
parameter values  where the map $T_{ab}$ has at least one
bounded orbit, which correspond to  $l_\infty$-eigenfunctions of
this difference operator. The paper shows 
that for  transcendental $\mu$ the
set $Spec_{\infty}[\mu]$ of energy values $E$ having a 
bounded solution is a Cantor set. Numerical
simulations suggest the possibility that these Cantor sets have positive
(one-dimensional) measure for all real values of $\mu$.

\vspace*{1.5\baselineskip}

\setlength{\baselineskip}{1.0\baselineskip}

%
%
%
%

\section{Introduction}
\hsp
As in parts I
and II, we  study the behavior under iteration of the two parameter 
family of piecewise-linear
homeomorphisms of $\RR^2$ given by
\beql{eq101}
T_{ab} (x,y) =
\left\{
\begin{array}{ccc}
(ax-y,x) & \mbox{if} & x \ge 0 , \\
(bx-y,x) & \mbox{if} & x < 0 .
\end{array}
\right.
\eeq
The parameter space is $(a,b) \in \RR^2$.
This map can be written
\beql{eq102}
T_{ab} (x,y) = \left[
\begin{array}{cc}
F_{ab} (x) & -1 \\ 1 & 0
\end{array}
\right]
\left[\begin{array}{c} x\\y\end{array}\right]
 ~,
\eeq
in which
\beql{eq103}
F_{ab} (x) = \left\{
\begin{array}{ccc}
a & \mbox{if} & x \ge 0 , \\
b & \mbox{if} & x < 0.
\end{array}
\right.
\eeq
The formula \eqn{eq102} shows that $T_{ab} (x,y)$ is a homeomorphism, 
since
\begin{equation}\label{eq104}
T_{ab}^{-1} (x,y) = 
 \left[ \begin{array}{cc}
F_{ab} (y) & -1 \\ 1 & 0
\end{array}\right]^{-1}
\left[\begin{array}{c} x\\y\end{array}\right]
=  \left[ \begin{array}{cc}
0 & 1 \\ -1 & F_{ab} (y)
\end{array}\right]
\left[\begin{array}{c} x\\y\end{array}\right]
\end{equation}
and it preserves the area form $d \omega = dx \wedge dy$.
In part I we observed that iterating this map 
encodes the solutions of the
second-order nonlinear recurrence
\beql{eq105}
x_{n+2} = \mu | x_{n+1} | + \nu x_{n+1} - x_n ~
\eeq
via
\beql{eq106}
T_{ab} (x_{n+1},x_n ) = (x_{n+2}, x_{n+1} )
\eeq
in which
\beql{eq107}
\mu= \frac{1}{2}(a - b), \qquad \nu:=\frac{1}{2}(a + b). 
\eeq
or equivalently 
$$
a = \nu + \mu , \quad b = \nu - \mu .
$$
This recurrence can be rewritten as a
one-dimensional discrete nonlinear  difference equation
of Schr\"{o}dinger type
\beql{eq108}
-x_{n+2} + 2x_{n+1} - x_n + V_\mu (x_{n+1} ) x_{n+1} = E x_{n+1}, ~
\eeq
where the potential $V_\mu (x)$ is given by
\beql{eq109}
V_\mu (x) := \left\{
\begin{array}{ccc}
\mu & \mbox{if} & x \ge 0 , \\
-\mu & \mbox{if} & x < 0.
\end{array}
\right.
\eeq
and the  energy value $E$ is given by
\beql{eq110}
E := 2- \nu ~.
\eeq
Holding the potential $V_\mu$ fixed and letting the parameter $\nu$ vary 
amounts to studying the set of solutions for
all energy values  $E = 2 - \nu$.

The values of $E$ corresponding to bounded orbits are
analogous to  the $l_\infty$-spectrum of
discrete linear Schr\"{o}dinger operators 
\beql{eq110a}
\Phi_V(x)_{n+1} = -x_{n+2} + 2x_{n+1} - x_n + V (n+1 ) x_{n+1}. 
\eeq
in which linearity is reflected in the 
potential $V(n)$ depending only on $n$ and not on $x_n$.
The model \eqn{eq110a}
is often called the ``tight-binding''
approximation to the Schr\"{o}dinger operator on the line. 
For a bounded potential $\Phi_V$ it gives
a  well-defined bounded operator 
on all sequence spaces $l_p(\ZZ)$. 
The eigenvalue equation is 
$$
\Phi(x)_{n+1} = -x_{n+2} + 2x_{n+1} - x_n + V (n+1 ) x_{n+1} = Ex_{n+1}.
$$
In  this context one is interested in characterizing 
the values of $E$ which have an orbit of the following types.

(1) {\em Extended state.} The orbit
$ \{ x_n : - \infty \le n \le \infty \}$ is a bounded orbit,
i.e. lies in $l_{\infty}(\ZZ)$.

(2) {\em Localized state.}
The orbit $ \{ x_n : - \infty \le n \le \infty \}$ lies
in $l_2(\ZZ)$.

\noindent The energy values $E$ for which there exists  an orbit of type (1)
comprise the {\em $l_\infty$-spectrum}
$$
Spec_{\infty}[\mu] := \{E= 2-\nu: (\mu, \nu) \in \Omega_{bdd} \}.
$$
A weaker version of (2) is the topological property that 
the orbit $\{ x_n : - \infty \le n \le \infty \}$
satisfies 
\beql{eq111}
\lim_{n \to \pm \infty} x_n =0 ~;
\eeq
we call such an orbit {\em weakly localized}.

Much work on the discrete linear  Schr\"{o}dinger operator was motivated
by the observation of Hofstadter \cite{Ho76} in 1976
that for  a  quasiperiodic linear potential
$V(n) = \lambda \cos (2 \pi \alpha n )$ 
with $\lambda=2$  
there is an $l_\infty$  eigenvalue structure $\Sigma_{\lambda, \alpha}$ 
which when (numerically) plotted for
variable  $\alpha$ appeared to 
form a fractal ``butterfly''.  More precisely, 
 for fixed irrational $\alpha$
the $l_\infty$-spectrum appeared to be a Cantor set of measure zero.
Hofstadter's model has been  much studied, and the
Hofstadter ``butterfly'' has been justified to some extent.
This has been done  
particularly in the context of the $l^2$-spectrun
of the almost Mathieu equation,
see Jitomirskaya \cite{Ji99} and  Puig  \cite{Pu04}
for recent results. 
Further references are given at the end of the introduction.

Here we  study analogous questions for the 
for the nonlinear 
difference operator of Schr\"{o}dinger type \eqn{eq108}.
We  obtain rigorous results about the $l_\infty$ spectrum
by exploiting the piecewise-linear structure of the maps.
The object of this paper is 
to determine structural properties of the set 
$$
\Omega_{SB} := \{ (\mu, \nu):~T_{\mu\nu} ~\mbox{has at least one nonzero 
bounded orbit} \},
$$
where  $SB$=''semi-bounded''.
We also obtain 
results concerning the set of parameter values with all orbits
bounded, which we denote 
$$
\Omega_{B} :=\{ (\mu, \nu):~T_{\mu\nu} ~\mbox{has all nonzero orbits
bounded} \}.
$$ 
Clearly $\Omega_{B} \subseteq \Omega_{SB},$ 
and Theorem~\ref{th26}(1) of this paper implies that
the inclusion is strict.
The set $\Omega_{SB}$
is an analogue of the Hofstadter ``butterfly'' set
in our context.
We prove
that $\Omega_{SB}$ is a closed set, 
and present evidence 
for the following conjecture. \\

\noindent{\bf Conjecture B.}
{\em The set $\Omega_{SB}$ has positive two-dimensional Lebesgue measure.
Furthermore, for each real $\mu$ the $l_\infty$-spectrum
$$ 
Spec_{\infty}[\mu] := \{ E=2 - \nu:~  (\mu,\nu) \in \Omega_{SB}   \}
$$
has positive one-dimensional Lebesgue measure.} \\

\noindent
The main result of the paper is to prove that  the set 
$Spec_{\infty}[\mu]$
is a Cantor set (totally disconnected perfect set) for
all parameter values $\mu$ outside  an exceptional set $\sE$
consisting entirely of  algebraic
numbers (Theorem~\ref{th29}). Thus
these sets $Spec_{\infty}[\mu]$ exhibit a property
ascribed to the Hofstadter ``butterfly''.
The value $\mu=0$ is exceptional, and 
$ Spec_{\infty}[0]$ is the entire interval  $ [0,4]$.
As far as we know $\mu=0$ might be the only value in the exceptional set;
if so the set $\Omega_{SB}$ for $\mu > 0$ would have the
structure  $(Cantor~ set)~\times~(line)$.
We  present numerical evidence
suggesting that for all values of $\mu$
(including the exceptional values) 
the set $ Spec_{\infty}[\mu]$
has positive one-dimensional Lebesgue measure.
If true, this  would contrast with the Hofstadter ``butterfly''.

The main difficulty in the proof of Theorem~\ref{th29}
is to show that for those $\mu$ that are not algebraic numbers
the set $Spec_{\infty}[\mu]$ is totally disconnected.
This reduces to showing that there are a dense set of
rationals in $[0, 1/2]$  that have 
nondegenerate  rotation intervals.
Establishing this requires nontrivial argument,
for part I showed
that for all non-algebraic $\mu$
 infinitely many rational rotation numbers always have degenerate
rotation intervals,  namely $r= \frac{2n-1}{4n}$
for $n \ge 2$. The main part of the proof of Theorem~\ref{th29}
shows that for each transcendental number $\mu$  all
rational rotation numbers $\frac{k}{p}$ in $(0, \frac{1}{2})$ 
with a  prime denominator $p$ have 
nondegenerate rotation intervals (Theorem~\ref{th82}).
The arguments used have a number-theoretic flavor.

In \S2 we state the main
results, and establish them in \S3-\S4.
In \S5 we make some  concluding remarks and 
formulate open questions.

Prior work on these maps includes Herman \cite[VIII]{He86},
and Beardon, Bullett and Rippon \cite{BBR95}. Both these
works studied, among other things, 
the set $\Omega_{B}$ where all orbits are
bounded. Their results were discussed in parts I and II,
and we remark on them later in this paper.

The results of this paper may be compared  with various
results concerning linear difference Schr\"{o}dinger operators. 
For general references on spectra of linear difference Schr\"{o}dinger
operators see Bougerol and Lacroix~\cite{BL85}
and Pasteur and Pigotin~\cite{PF92}. For discussion of
existence or nonexistence of Cantor set spectra for various 
potentials see  Bellissard~\cite{Be92},
Fr\"{o}hlich et al \cite{FSW90} and the earlier references for
the almost Mathieu equation.
For Schr\"{o}dinger operator potentials taking
finitely many values, see Kotani~\cite{Ko90},
Sutherland and Kohmoto~\cite{SK87} and S\"{u}t\"{o}\cite{Su95}.

\noindent \paragraph{Notation.} We write $\bv=(\bv_x, \bv_y) \in \RR^2$,
to be viewed as a column vector.
 An interval $[\bv_1, \bv_2)$ of the 
unit circle, or a  corresponding sector $\RR^{+}[\bv_1, \bv_2)$
of the plane $\RR^2$, is specified by rotating  counterclockwise from
$\bv_1$ to $\bv_2$.
We let $\mbox{Meas}_d (S)$ denote the $d$-dimensional Lebesgue measure
of a set $S$, for $d= 1, 2$.

\noindent \paragraph{Acknowledgment.} We did much of 
the work on this paper while employed at AT\&T Labs-Research,
whom we thank for support. We thank
Jim Reeds for suggesting a proof method for Lemma~\ref{le82}.
We thank M. Kontsevich for bringing the work of
Bedford, Bullett and Rippon \cite{BBR95}
to our attention.

%
%
%
%
\section{Summary of Results}
\hsp
The parameter space of the map can be taken to be either
$(a, b)$ or $(\mu,\nu)$, as these are equivalent by 
\beql{eq201}
\mu = \frac{1}{2}( a - b), \qquad \nu= \frac{1}{2} (a + b).
\eeq
Both coordinate systems have their advantages, and we
 write the map \eqn{eq101} as $T_{ab}$, $T_{\mu \nu}$ accordingly.
We will also use the mixed parameter space $(a, \mu)$.
It is convenient to represent the action of $T_{ab}$, 
acting on row vectors $\bv_n = (x_{n+1}, x_n )$ as
\beql{eq202}
T_n (\bv_0) = (x_{n+1}, x_n)
 = M_n (\bv_0) (x_1,  x_0)
\eeq
in which
\beql{eq203}
M_n (\bv_0) = \prod_{i=1}^n \left[
\begin{array}{cc}
F_{ab} (x_i) & -1 \\ 1 & 0
\end{array}
\right] :=
\left[
\begin{array}{cc}
F_{ab} (x_n) & -1 \\ 1 & 0
\end{array} 
\right] \cdots
\left[
\begin{array}{cc}
F_{ab} (x_2 ) & -1 \\ 1 & 0 
\end{array}
\right] 
\left[
\begin{array}{cc}
F_{ab} (x_1 ) & -1 \\
1 & 0
\end{array}
\right]\,.
\eeq

Conjugation by the involution $J_0: (x,y) \to (-x, -y)$ gives
\beql{eq204}
T_{ba} = J_0^{-1} \circ T_{ab} \circ J_0 ~.
\eeq
Thus, in studying dynamics, without loss of generality we can
restrict to the closed half-space
$\{(a,b) : a \ge b \}~$ of the $(a,b)$ parameter space.
This corresponds to the region $\{(\mu, \nu): \mu \ge 0 \}$
of the $(\mu. \nu)$ parameter space, with $T_{\mu \nu}$ conjugate
to $T_{-\mu, \nu}$.

The associated rotation map $S_{\mu\nu}: S^1 \to S^1$ is given
by 
$$ S_{\mu\nu}(e^{i\theta}):= 
\frac{T_{\mu\nu}(e^{i\theta})}{|T_{\mu\nu}(e^{i\theta})|}.
$$
It has a well-defined rotation number $r(S_{\mu\nu})$, which
is a counterclockwise rotation, and was shown in part I to
always lie in the closed interval $[0, 1/2].$

In \S3 we establish properties of the sets $\Omega_{SB}$ and $\Omega_{B}$.
%
%
%

\begin{theorem}\label{th26}
The set $\Omega_{SB}$ is a closed set.
It consists of all parameter values $(\mu, \nu)$ for which 
the associated values $(a,b)$ satisfy one of the 
conditions below.

(1) $r(S_{ab} )$ is rational, and $T_{ab}$ has a periodic orbit.

(2) $r(S_{ab} )$ is irrational.
\end{theorem}

\noindent
Case (1) was already established in part I, Theorem 2.4(i) and (iii).
To handle case (2), we study the $(\mu, \nu )$-parameterization
 for constant $\mu$, and show the following facts.
(Theorem \ref{th62}.)

(a) The set of values $\nu$ with $r(S_{\mu \nu} ) = r$ with 
$0 < r= \frac{p}{q} < \frac{1}{2}$ is either a single point $\bv$ or 
a closed interval
$[\bv_1, \bv_2]$.
If it is a point then $T_{\mu \nu}$ is periodic.
If it is an interval $[\bv^- (r), \bv^+ (r)]$ then no $T_{\mu \nu}$ 
is periodic, and the only $T_{\mu \nu}$ with a periodic orbit 
are the endpoints
$\bv = \bv^- (r)$ and $\bv^+ (r)$.

(b) The set of values $\nu$ with $r(S_{\mu \nu} ) = r$ 
irrational with $0 < r < \frac{1}{2}$ is a single point $\bv$.

We then obtain a bounded orbit in case (b) by a limiting procedure
using bounded orbits in case (a) with suitable
rational rotations $\frac{p_n}{q_n}$
approaching $r$.  The endpoints of   rotation intervals
of positive length in (a) give points 
in $\Omega_{SB}$ that are not in $\Omega_{B}$.

We obtain the following information on the location of $\Omega_{SB}$ 
viewed in the
$(\mu. a)$ parameter space, using results from \S3 of part I.

%
%
%
\begin{theorem}\label{th27}
The set of values for which $T_{\mu\nu}$ 
has a bounded orbit (i.e. $(\mu,\nu) \in \Omega_{SB}$),
in the range $\mu \ge 0$ 
lies inside the 
cylinder $-2 \le a \le 2 ~,$ where $a= \nu + \mu$.
Furthermore:

(1) For $-2 \le a < 0$,  one has
$$
0 \le \mu \le \frac{2}{|a|}- \frac{|a|}{2}.
$$

(2) For $0 \le a \le 2$, set $a = 2 \cos \theta$.
Then for 
for each $n \ge 2$ and $\frac{\pi}{n+1} \le \theta < \frac{\pi}{n}$,
there holds 
$$
0 \le   \mu   \le \cos \theta - 
\frac{\sin n\theta + \sin \theta}{\sin (n+1)\theta}.
$$

\end{theorem}

\begin{figure}[htb]
\centerline{\psfig{file=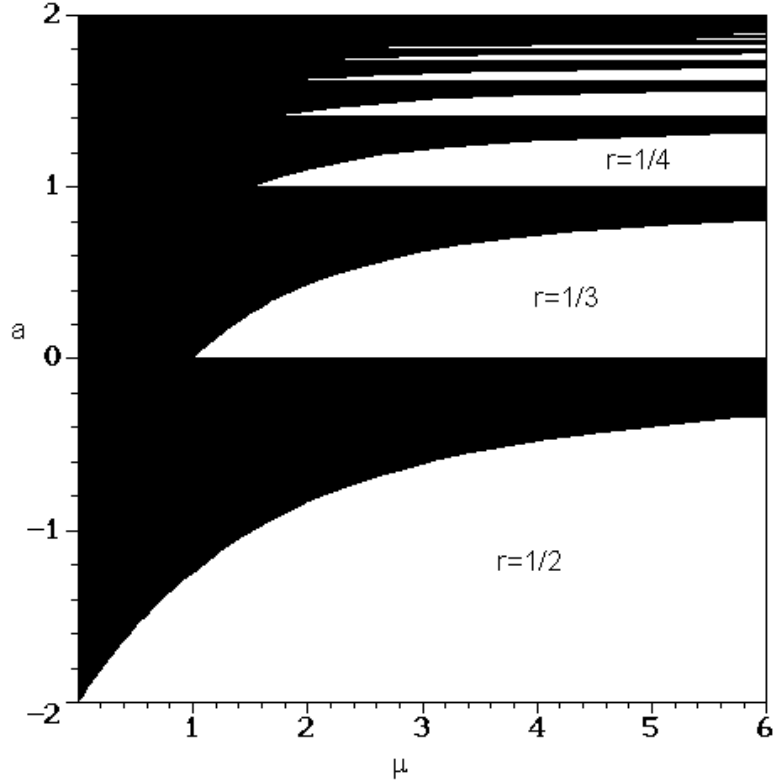,width=4.5in}}
\caption{Region $\sR$ containing $\Omega_{SB}$ in $(a,\mu)$-plane.}
\label{fg2.1}
\end{figure}

\noindent
A weaker bound for $\Omega_{B}$ was previously obtained by
Beardon, Bullett and Rippon \cite[Theorem 1.1(iii)]{BBR95}.
In Figure \ref{fg2.1} below
we plot the parameter region allowed in Theorem \ref{th27} 
in the $(\mu,a)$-plane, with $\mu = \frac{1}{2} (a-b)$, for 
$0 \le \mu \le 6$. 
%
%
The region $\sR$ of Theorem \ref{th27} becomes unbounded along the lines
$$
a= a_n = 2 \cos \frac{\pi}{n} , \quad n \ge 2 \,,
$$
which are $a=0,1, \sqrt{2}, 2 \cos \frac{\pi}{5}$, $\sqrt{3} \ldots$.

The excluded regions in the strip $-2 \le a \le 2$ touching 
$a_n$ for $n \ge 1$ are regions where $T_{\mu\nu }$ has rational 
rotation number $r_n = \frac{1}{n}$.
At the end of \S3 we present evidence that the set $\Omega_{SB}$
has positive Lebesgue measure, and that each
$Spec_{\infty}[\mu]$ has positive one-dimensional Lebesgue measure.
The bounds of Theorem ~\ref{th27} imply that
$$
\lim_{\mu \to \infty} ~\mbox{Meas}_1( Spec_{\infty}[\mu]) = 0.
$$
The regions associated to the rotation numbers $r_n = \frac{1}{n}$
asymptotically remove all the area in $-2 \le a \le 2$.

For comparison with Theorem~\ref{th26}, 
we include  the following result 
characterizing the set $\Omega_{B}$ of all parameter
values for which all
orbits are bounded, given in 
 Beardon, Bullett and Rippon \cite{BBR95}, based on 
results of Herman \cite[VIII2.4]{He86}.

%
%
%
\begin{theorem}\label{th28}
The set $\Omega_{B}$ 
consists of those 
parameter values $(\mu, \nu)$
for whose  associated  $(a,b)$ the map 
$T_{ab}$ is conjugate to a rotation of
the plane. This occurs if and only if one of the
following occur.

(1) $r(S_{ab} )$ is rational and $T_{ab}$ is of finite order, 
i.e. $T_{ab}^{(k)} =I$ for some finite $k$.

(2) $r(S_{ab} )$ is irrational, and $T_{ab}$ contains an invariant circle.
\end{theorem}

Examples given in part I and II show that cases (1) and (2) both occur. 
The set of parameter values for which case (1) occurs has
Hausdorff dimension $1$. 

In \S4  we study bounded solutions for the  one-parameter ``eigenvalue''
families in which  $\mu = \mbox{constant}$. 
Recall that a {\em Cantor set} in $\RR$ is a perfect\footnote{A closed set 
$E$ is {\em perfect} if each $x \in E$ is a limit point 
of a infinite sequence of distinct points of $E$.} totally disconnected set.
The following theorem is the main result of the paper.
%
%
%
\begin{theorem}\label{th29}
For $\mu_0 \in \RR$ the set
$$ 
Spec_{\infty}[\mu] := \{ E=2 - \nu:~  (\mu,\nu) \in \Omega_{SB}   \}
$$
is a Cantor set if $\mu$ is not in a countable exceptional set 
$\sE$ that consists entirely of algebraic numbers.
\end{theorem}

The set $\sE$ contains $\mu=0$, where $ Spec_{\infty}[0]  = [0, 4]$.
As far as we know, this might be the only point in $\sE$.
If so, then  $\Omega_{SB}$
would have  the structure $(Cantor~ set) \times (half{-}line)$ in
the $(\mu,a)$-parameter space, in the region $-2 \le a \le 2$ and
$\mu > 0$. 
As indicated earlier, 
numerical evidence supports the assertion that  
$Spec_{\infty} [ \mu ]$ has positive one-dimensional Lebesgue measure
for all real values of $\mu$.

%
%
%
%
\section{Bounded Orbits}
\hsp
We consider the set $\Omega_{SB}$ of parameter values 
having at least one bounded orbit.

\setcounter{figure}{0}
%
%
%

\begin{theorem}\label{th61}
The set $\Omega_{SB}$ is a closed set.
Each  $(\mu, \nu) \in \Omega_{SB}$ 
the associated $T_{ab}$ possesses a bounded orbit 
$\sO (\bv_0) = \{\bv_n : k \in \ZZ \}$ such that
\beql{eq61}
\| \bv_0 \| = \sup_{n \in \ZZ} ~
\|T_{ab}^{(n)} (\bv_0 ) \| ~.
\eeq
\end{theorem}

\paragraph{Proof.}
Let $(a_k,b_k )$ be associated values to
parameters $(\mu_k, \nu_k) \in \Omega_{SB}$
with $(a_k, b_k ) \to (a,b)$.
We can choose a bounded orbit
$\sO (\bv_0^{(k)} ) = \{ \bv_n^{(k)}: n \in \ZZ \}$
of $T_{a_k b_k}$ with
$$\sup_{n \in \ZZ} \| \bv_n^{(k)} \| = 1\,,$$
and we may also suppose that
$$1 \ge \| \bv_0^{(k)} \| \ge 1 - \frac{1}{k}\,,$$
by shifting the orbit appropriately.
By compactness of the unit ball in $\RR^2$ we may extract
a subsequence so that
$$\bv_0^{(k_j)} \to \bv_0 ~,$$
and necessarily $\| \bv_0 \| =1$.
Set $\bv_m = T_{ab}^{(m)} (\bv_0)$.
Then, for each fixed $m \in \ZZ$,
$$\bv_m^{(k_j)} \to T_{ab}^{(m)} (\bv_0) ~,
$$
as $k_j \to \infty$, because $T_{ab} (\cdot )$ depends continuously
on the parameters
$(a,b) \in \RR^2$.
Since
\linebreak
$\| \bv_m^{(k_j)} \| \le 1$, this yields
$$\| \bv_m \| = \| T_{ab}^{(m)} (\bv_0 ) \| \le 1 ~.$$
Thus $T_{ab}$ has a bounded orbit, so $(a,b) \in \Omega_{SB}$, 
and furthermore this orbit attains its supremum
$$\| \bv_0 \| = \sup_{m \in \ZZ} \| \bv_m \| =1 ~.$$
The argument applies to any $(a,b)$ associated
to some $(\mu, \nu) \in \Omega_{SB}$ by taking all
$(a_k, b_k) = (a,b)$, to give \eqn{eq61}.~~~$\bsq$

To prove Theorem \ref{th26} we establish the following preliminary result,
 using $(\mu,\nu)$-parameters.

%
%
%

\begin{theorem}\label{th62}
Let $\mu \in \RR$ be fixed and let $\nu$ vary over $-\infty < \nu < \infty$.

(1) Let $0 < r <  \frac{1}{2}$ be rational. Then the set $I_{\mu}(r)$ 
of values $\nu$ such that $r(S_{\mu \nu} ) = r$   
is either  a point
$\nu^{\pm}(r)$ or an interval $[\nu^- (r) , \nu^+ (r) ]$.
In the point case $T_{\mu \nu^{\pm}}$ is of finite order.
In the interval case $T_{\mu \nu}$ is never of finite order, 
 and $T_{\mu \nu}$ for $\nu \in [\nu^- (r), \nu^+ (r) ]$ contains 
a periodic orbit if and only if $\nu$ is one of the endpoints of the interval
$\nu = \nu^- (r)$ or $\nu^+ (r)$.

(2) Let $0 < r <  \frac{1}{2}$ be irrational.
Then the  set  $I_{\mu}(r)$   of values $\nu$ 
such that $r(S_{\mu \nu} ) = r$ 
consists of a point $\nu^\pm (r)$.

(3) The set $I_{\mu}(0)$ of values $\nu$ with $r(S_{\mu \nu} ) = 0$ is a 
half-infinite interval $[\nu^- (0), + \infty )$.
In this interval $T_{\mu \nu }$ is never periodic, and 
$T_{\mu \nu}$ has a periodic orbit only for $\nu = \nu^- (0)$.

(4) The set  $I_{\mu}(1/2)$  of values 
$\nu$ with $r(S_{\mu \nu}) = \frac{1}{2}$ 
is a half-infinite interval
$(-\infty, \nu^+ (\frac{1}{2} )]$.
In this interval $T_{\mu \nu}$ is never periodic, and 
$T_{\mu \nu}$ has a periodic orbit only for $\nu = \nu^+ (\frac{1}{2} )$.
\end{theorem}

\paragraph{Proof.}
The continuity and nonincreasing properties of $S_{\mu \nu}$ 
in Theorem 2.2 of part I
imply that for $0 \le r \le 1/2$ the set
$\{ \nu : r (S_{\mu \nu} ) =r \}$ is either a point or an interval.
We first consider cases (1), (3) and (4).
Lemma 4.1 of part I
implies that for rational 
$r(S_{\mu \nu} ) = \frac{p}{q}$, 
the point case can occur if and only if $S_{\mu \nu}$ is periodic, 
and by Theorem 2.5 (iii) of part I,
this occurs if and only if
$T_{\mu \nu}$ is of finite order.
Thus in the interval case
$I(r) =[\nu^- (r) , \nu^+ (r)]$, the map $T_{\mu \nu}$ is never periodic, 
and the classification of Theorem 2.4 of part I
showed that the circle map 
$S_{\mu \nu}$ then has either one or two periodic orbits.
The proof of Theorem 2.4 of part I
showed that in the interior of 
the rotation interval $(\nu^- (r), \nu^+ (r))$ the 
circle map $S_{\mu \nu}$ has two periodic orbits whose points 
alternate around the circle
$0 \le \theta \le 2 \pi$.
At the endpoint values $\nu^+ (r)$, $\nu^- (r)$ these 
coalesce into a single periodic orbit, which is case (i) of 
Theorem 2.4 of part I,
and $T_{\mu \nu}$ then has a periodic orbit.
This proves (1).
Cases (3) and (4) follow similarly.

Now consider case (2), where $r(S_{\mu \nu_0} ) =r$ is irrational.
We show first that for any $ \ep > 0$ there exists 
$\nu'$ in $[\nu_0 , \nu_0 + \ep ]$ such that $S_{\mu \nu'}$ 
has a periodic point, 
so that $r(S_{\mu \nu'} )$ is rational.
The proof of Lemma 4.1 of part I
showed that
$$ 
\left. \frac{\partial}{\partial \nu} 
S_{\mu \nu}^{(2)} (\theta )\right|_{\nu = \nu_0} < 0~.
$$
Thus, choosing $\nu_1=\nu_0 + \ep$ there is a positive constant 
$\ep_1$ such that the lifted map $\tilde{S}_{\mu \nu}^{(2)}$ to the line has
$$\tilde{S}_{\mu \nu_1}^{(2)} (\theta ) \le
S_{\mu \nu_0}^{(2)} (\theta ) - \ep_1, \quad
0 \le \theta \le 2 \pi
$$
since $S_{\mu \nu}^{(2)}$ is strictly monotone in $\theta$;
this gives
\beql{eq64}
\tilde{S}_{\mu\nu_1}^{(2)} (\theta ) \le 
\tilde{S}_{\mu \nu_0}^{(2)} ( \theta ) = \ep_1 , \quad\mbox{all} \quad
\theta \in \RR \,.
\eeq
Since $R(S_{\mu \nu_0} ) =r$ is irrational, 
one can find a positive integer $m$ such that
\beql{eq65}
0 < S_{\mu \nu_0}^{(m)} (\theta ) < \ep_1 \qquad (\bmod~2 \pi ) ~.
\eeq
Then as $\nu$ increases continuously in $[\nu_0 , \nu_0 + \ep ]$,
the endpoint $S_{\mu \nu_0}^{(m)} (\theta )$ moves 
continuously and monotonically downward, 
there exists some $v' \in [\nu_0 , \nu_1 ]$ with
\beql{eq66}
S_{\mu\nu}^{(m)} (\theta ) =0 \qquad ( \bmod~2\pi ) ~,
\eeq
by \eqn{eq64} the movement of 
$S_{\mu \nu}^{(m)} (\theta )$ as $\nu$ varies is 
counterclockwise by at least $\ep_1$.
Now \eqn{eq66} shows that $S_{\mu \nu'}$ has a 
periodic point of period $m$, hence $R(S_{\mu \nu'} )$ 
is rational, as claimed.
It follows that $r(S_{\mu \nu_0}) > r(S_{\mu \nu '})$ 
and $\nu_0 < \nu ' \le \nu_0 + \ep$.

Next, a similar argument shows that for each $\ep > 0$ there exists
$\nu ''$ in $[\nu_0 - \ep , \nu_0 ]$ such that $S_{\mu \nu''}$ 
has a periodic point.
Thus $r(S_{\mu \nu_0}) < r(S_{\mu \nu''} )$ with 
$\nu_0 - \ep \le \nu'' \le \nu_0$.

We conclude that the rotation interval $I(r)$ for 
the irrational rotation value $r$
(for fixed $\mu$) is contained in 
$[\nu_0 - \ep , \nu_0 + \ep ]$ for all $\ep > 0$, 
hence $I(r) = \{\nu_0 \}$ is a point, which we label
$\nu^{\pm} (r)$.~~~$\bsq$ \\

%
%
%

\begin{theorem}~\label{th61b}
If the rotation number $r(S_{ab})$ is irrational,
then $T_{ab}$ has a bounded orbit. 
\end{theorem}

\paragraph{Proof.}
We switch to $(\mu, \nu )$ parametrization, and suppose that
$r(S_{\mu \nu_0}) = r$ is irrational,
with $\mu = \frac{1}{2} (a-b)$,
$\nu_0= \frac{1}{2} (a+b)$;
necessarily $0 < r < \frac{1}{2}$.
We consider $T_{\mu \nu}$ for variable $\nu$, and show we can 
pick a sequence $\{ \nu_k : k \ge 1 \}$ with $\nu_k \to \nu_0$ such that
\begin{quote}
$(\ast )$ $T_{\mu\nu_k}$ has rational rotation number and a 
periodic orbit $\{ \bv_n^{(k)} : n \in \ZZ \}$, normalized with
$\| \bv_0^{(k)} \| = \max_{n \in \ZZ} \| \bv_n^{(k)} \| =1$.
\end{quote}

\noindent
To do this we pick a series of rational approximations
$\{ \frac{p_k}{q_k} \}$ approaching $r$ monotonically from above,
and choose $\nu_k = \nu^- (\frac{p_k}{q_k})$, for the rotation interval
$I(r) = \left[ \nu^- (\frac{p_k}{q_k} ), \nu^+ (\frac{p_k}{q_k} ) \right]$.
By Theorem \ref{th62} (1) $T_{\mu \nu_k}$ has a periodic orbit, 
which we can normalize by scaling to satisfy $(\ast )$.
By compactness, we can extract a subsequence
$\{ \bv_0^{(k)} \}$ such that $\bv_0^{(k_j )} \to \bv_0$
as $j \to \infty$.
Then $\| \bv_0 \| =1$ and, just as in the proof of
Theorem \ref{th61},
$\sO(\bv_0) = \{ \bv_n = T_{\mu \nu_0}^{(n)} (\bv_0) : n \in \ZZ \}$ 
is a bounded orbit of $T_{\mu \nu_0}$, with
$$\| \bv_0 \| = \max_{n \in \ZZ} \| \bv_n \| =1 ~,$$
as required.~~~$\bsq$ \\

We now deduce Theorems  \ref{th26} and  \ref{th27}.

\vspace*{+.1in}
\noindent
%
%
%

\paragraph{Proof of Theorem \ref{th26}.}

(1) Suppose that $r(S_{ab} )$ is rational.
Theorem 2.4 (i) and (iii) of part I classify when $T_{ab}$ 
has a bounded orbit.
All such orbits are periodic orbits of $T_{ab}$.

(2) The case  that $r(S_{ab} )$ is irrational
is covered by  Theorem~\ref{th61b}. 
~~~$\bsq$

%
%
%

\paragraph{Proof of Theorem \ref{th27}.}
The excluded regions in the theorem, for $-2 \le a \le 0$ correspond
to rotation number $\frac{1}{2}$ and for $a= 2\cos \theta \ge 0$ and
$ \frac{\pi}{n+1} < \theta < \frac{\pi}{n}$ for $n \ge 2$ they
correspond to rotation number $\frac{1}{n+1}$.
The proof of Theorem 3.3 of part I  shows that for $a \ge b$
(that is, $\mu \ge 0$) 
and parameter values outside the stated range $\sR$, 
$S_{\mu\nu}$ has an
associated matrix $M$ with $Tr(M) > 2$, and has
two periodic orbits, corresponding to the two real eigenvectors of $M$.
Case (2) of  Theorem 2.4 then yields that every orbit of $T_{\mu\nu}$
is unbounded.~~~$\bsq$ \\

\begin{figure}[htb]
\centerline{\psfig{file=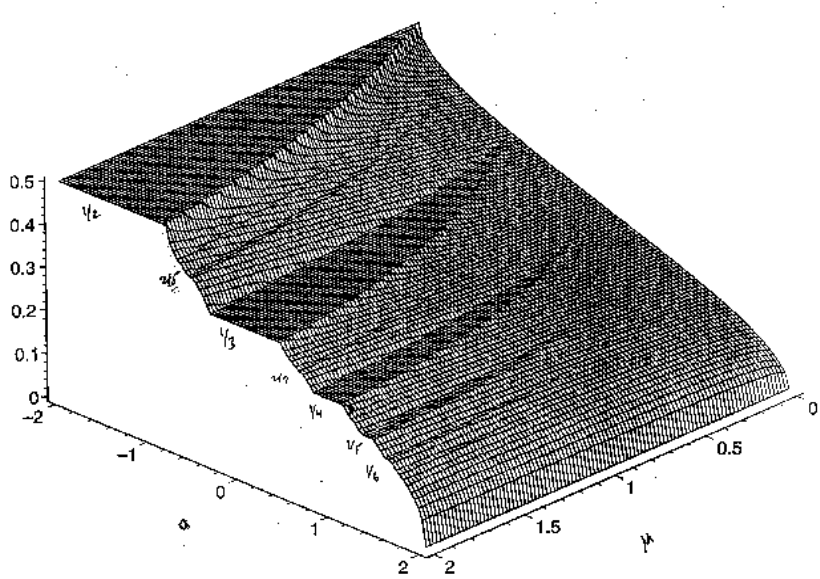,width=4.5in}}
\caption{Rotation number in $(\mu , a)$-plane, $0 \le \mu \le 2$, 
$-2 \le a \le 2$.}
\label{fg6.1}
\end{figure}

The allowed region $\sR$ in Theorem \ref{th27} containing $\Omega_{SB}$ was
pictured in \S2, in Figure \ref{fg2.1}.
It can be compared with Figure \ref{fg6.1} below
where we numerically plot 
the rotation number $r(T_{ab} )$ in the $(\mu, a )$-parameter plane, 
where $\mu \ge 0$
and $-2\le a \le 2$, in the range $0 \le \mu \le 2$.
%
%

\noindent The rotation intervals for 
$r = \frac{1}{2}$, $\frac{2}{5}$, $\frac{1}{3}$,
$\frac{2}{7}$, $\frac{1}{4}$, $\frac{1}{5}$ and $\frac{1}{6}$ 
are visible in Figure \ref{fg6.1}. 
We note that the rotation interval with $r = 1/3$ pictured above
is larger than that excluded from the region $\sR$,
for Theorem~\ref{th26} excludes no values 
for $0 \le \mu \le 1$ with $0 \le \alpha < 1,$
corresponding to $r = \frac{1}{3}$.

We now describe numerical evidence suggesting that $\Omega_{SB}$ has 
positive Lebesgue measure, in the region $0 \le \mu \le 1$.
The question is equivalent to:
Does the (open) set consisting of the interiors of all rational 
rotation intervals have full Lebesgue measure in the region
$-2 \le a \le 2$, $0 \le \mu \le 1$?
We suggest it does not.
Consider, for fixed $\mu$, the amount of area in $\{-2 \le a \le 2\}$
left uncovered by the interiors of rational rotation intervals.
At $\mu =0$ the full interval $[-2, 2]$ of measure 4 is left uncovered.
(For $\mu =0$ all rotation intervals have width 0, 
by Example 3.1 of part I.)
We observe empirically in Figure~\ref{fg6.1}  that the width 
of rational rotation numbers $r=\frac{1}{n}$ appears to be
nondecreasing in $\mu$. We know from Example 3.3 of part I
that for certain values, such as $r = \frac{3}{8}$, 
the rotation interval has width zero for $0 \le \mu \le \infty$.
There are rational rotation numbers (e.g. $r = \frac{2}{9}$)
where the length of the rotation interval as a
function of $\mu$ is sometimes increasing
and sometimes decreasing.
For small $\mu$ we found 
empirically that 
the total uncovered length appears nonincreasing as a function
of $\mu$. If this holds up to $\mu=1$, then in the range
$0 \le \mu \le 1$, the uncovered length would be minimized at $\mu =1$.
For $\mu =1$ we numerically estimated the length in 
$\{-2 \le a \le 2\}$ left uncovered by rotation intervals 
$r = \frac{p}{q}$ with denominator $q \le N$ and obtained the 
data in the following table. 
\begin{table}[htb]
\begin{center}
\begin{tabular}{|l|l|} \hline
\multicolumn{1}{|c|}{Denominator} & \multicolumn{1}{c|}{Area} \\ \hline
N=1 & 4.00000000 \\
N=2 & 3.23606798 \\
N=4 & 2.45423817 \\
N=8 & 2.17951560 \\
N=16 & 2.09869137 \\
N=32 & 2.05071812 \\
N=64 & 2.02839432 \\
N=128 & 2.01690968 \\
N=256 & 2.01094576 \\ \hline
\end{tabular}
\end{center}
\caption{Uncovered length: case  $\mu=1$.}
\end{table}

\noindent This data suggests that in the limit as the denominator cutoff
$N \to \infty$, for $\mu=1$
there is an uncovered set of positive Lebesgue measure,
and this measure is approximately 2.
Under the monotonicity hypothesis above,
the two-dimensional Lebesgue measure of 
$\Omega_{SB}$ in the region $0 \le \mu \le 1$ in the $(\mu , a )$ 
plane would lie between $2$ and $4$.

The numerical  evidence above led us to formulate Conjecture B.
As additional evidence in favor of Conjecture B, we 
numerically located orbits that appear to be invariant 
circles, described in \S5 of part I.
 If invariant circles were a measure zero phenomenon, we would not 
expect to numerically find such invariant circles. 
It appears plausible that $\Omega_{B}$ may also have positive 
(two-dimensional) Lebesgue measure.

We now prove Theorem~\ref{th28}, which amounts to
the equality
$\Omega_{R}= \Omega_{B}$, where $\Omega_{R}$
are the parameter values for which $T_{ab}$ is
conjugate to a rotation of the plane. We follow 
Beardon, Bullett and Rippon \cite{BBR95}.

%
%
%

\paragraph{Proof of Theorem \ref{th28}.}
Case (1) is covered by Theorem 2.4 of part I. Case (2) follows
from results of Herman \cite[VIII.2.4]{He86}.
For $S_{ab}$ having  irrational rotation number 
he establishes the dichotomy that either 
$T_{ab}$ has an invariant circle, or else it is
topologically transitive, i.e. it contains an
orbit that is dense in $\RR^2$. Such an orbit is unbounded.

If $T_{ab}$ is topologically conjugate to a rotation of 
the plane,  then all orbits are bounded. 
The fact that $T_{ab}$ is topologically conjugate to a rotation
in cases (1) and  (2)  follows 
also from results of Herman \cite[VIII.2.5]{He86},
which uses facts from Herman \cite{He79}.
~~~$\bsq$ \\

Finally we consider the possible
existence of localized orbits. We call an orbit 
$\{ \bx_n: n \in \ZZ\}$
{\em weakly localized} if
\beql{350}
\lim_{n \to \pm \infty} ||\bx_n|| =0.
\eeq
Theorem~\ref{th28} implies that if $T_{\mu\nu}$ had
a localized orbit, then necessarily $S_{\mu\nu}$ has irrational
rotation number, and $T_{\mu\nu}$ is not conjugate to a rotation.
It is not known whether any $T_{\mu\nu}$ with the 
latter property exists. Herman \cite{He86} observed that
any such map necessarily has a dense orbit. 
Herman \cite[VIII.2.4]{He86} conjectured and gave
evidence for the
existence of some element of the Froeschl\'{e} group
having this property. Here the {\it  Froeschl\'{e} group} is the
group of homeomorphisms of the plane generated by
$SL(2, \RR)$ and all maps
$$
G_{ab}= \left[ \begin{array}{cc}
1 & 0 \\
F_{ab}(x) & 1 
\end{array} \right],
$$
for real $a, b$, viewed as acting on column vectors.
Bedford, Bullett and Rippon 
\cite[pp. 673-674]{BBR95} discuss the possibility 
that this occurs for some  maps $T_{ab}$, and suggest it might even
occur somewhere on the parameter line  $\nu=1/2(a+b)=0$.

%
%
%
%
 
\section{Cantor Set Spectra: Constant $\mu$}  
\hsp
We consider the family $T_{\mu\nu}$ for fixed $\mu$ and
variable $\nu$. Recall that 
for these families, one is particularly interested in knowing for which
$E=2-\nu$ one has bounded orbits; these values of $E$ correspond
to eigenvalues of the nonlinear Schr\"{o}dinger system (1.8)
giving  extended states or localized states.  For 
the case $\mu=0$ where the maps $T_{\mu\nu}$ 
are linear  we have:
%
%
%

\begin{theorem}\label{th81}
For $\mu=0$, $T_{\mu\nu}$ has a nontrivial bounded orbit
if and only if $\nu\in[-2,2]$; all orbits are bounded if and only if
$\nu\in(-2,2)$.
\end{theorem}

\paragraph{Proof.}
If $\nu\not\in[-2,2]$, then, by Theorem \ref{th26}, $T_{\mu\nu}$ has
no bounded orbits.  If $\nu=\pm 2$, then
explicit calculation verifies that $S_{\mu\nu}$ has exactly 
one periodic orbit;
by Theorem 2.4 (i) of part I, 
then $T_{\mu\nu}$ has exactly one nontrivial bounded orbit.
The remaining case, $\nu\in(-2,2)$, was covered in Example 3.1 of part I
which showed it is conjugate to a rotation.  
In this case $T_{\mu\nu}$ either is periodic or else 
has an invariant circle (an
ellipse, in fact).~~~$\bsq$
\bigskip

For fixed $\mu \in \RR$ say that a rational rotation interval 
$I_\mu (r)$ in the $\nu$-variable is {\em nondegenerate} if
$I_\mu (r) = [\mu_\mu^- (r), \nu_\mu^+ (r)]$, has positive measure and is
{\em degenerate} if $I_\mu (r) = [ \nu_\mu^0 (r)]$ is a point.
We let 
$$
E(\frac{m}{n}) :=\{ \mu \ge 0: I_\mu (\frac{m}{n}) ~~\mbox{is degenerate}\}.
$$
denote the ``'exceptional'' set where degeneracy occurs.

To show that 
$Spec_{\infty}[\mu] := \{ E= 2 - \nu : 
T_{\mu \nu} ~\mbox{has a bounded orbit} \}$ 
is a Cantor set for a given  $\mu$, we must show that it is
totally disconnected; for this it is necessary and sufficient that there  are
a dense set of rational rotation numbers in $[0, 1/2]$ whose
rotation intervals at $\mu$ are nondegenerate. The
following result gives a criterion for nondegeneracy which
applies to rotation numbers having an odd prime denominator.

%
%
%
\begin{theorem}\label{th82}
Let $q \ge 3$ be prime. For each rational 
$r= \frac{j}{q}$ with $ 1 \le j \le \frac{q-1}{2}\}$ the set
 $E( \frac{j}{q})$  of all
$\mu \ge 0$ for which the rotation interval of rotation number
$\frac{j}{q}$ is degenerate is a 
finite set of algebraic numbers.
\end{theorem}

\noindent\paragraph{Remark.} Each set $E( \frac{j}{q} )$ includes $\mu =0$, 
and may include other points. For example,  $E(\frac{2}{7} )$ 
contains $\mu = \frac{1}{2}$. Theorem~\ref{th82}
does not generalize to hold for all $E(\frac{m}{n})$; 
indeed  Example 3.3 of part I 
shows that $E(\frac{1}{8})$ is the positive real line.\\

To prove this result we first establish some preliminary lemmas.
We shall describe iteration of the map $T_{ab}$
assuming the iterates have  a given symbol sequnece
$\sS=(S_1, S_2, ..., S_n)$, with $S_i = \pm 1$ instructing whether
$F_{ab}(x) = ax$ or $bx$. We keep track of the iterates for
a symbol sequence  
using $(a, \mu)$ coordinates, where $a$ and $\mu$ are viewed
as indeterminates; note that $b= a - 2 \mu$. Requiring 
periodicity of a particular orbit  then
gives polynomial  equations that $a$ and $\mu$ 
must satisfy. We first treat the case corresponding to  $\mu \ne 0$.

%
%
%
\begin{lemma}\label{le82}
Given a symbol sequence 
$\sS = (S_1, \ldots, S_n ) \in \{1, -1\}^n$
define bivariate polynomials
$p_{n}^\sS (a, \mu) \in \ZZ [a, \mu]$ by
\beql{eq85}
\left[
\begin{array}{c}
p_n^\sS (a, \mu ) \\ p_{n-1}^\sS (a,\mu ))\end{array}\right] :=
\left[ \begin{array}{cc}
F_{ab } (S_{n-1}) & -1 \\1 & 0 \end{array}\right] \cdots
\left[ \begin{array}{cc}
F_{ab} (S_1 ) & -1 \\ 1 & 0
\end{array}\right]
\left[\begin{array}{c}1\\0\end{array}\right]
,
\eeq
where 
\beql{eq86}
F_{ab} (S_i) = \left\{
\begin{array}{ccc}
a & \mbox{if} & S_i =1,  \\
a-2\mu & \mbox{if} & S_i =-1. 
\end{array}
\right.
\eeq
Now suppose $\mu \ne 0$ is a fixed real number,
 and view each $p_{n}^\sS (a, \mu)$
as a univariate polynomial in $\RR[a]$. 
If  $\sS, \sS'$ are two distinct symbol sequences, then
\beql{eq87}
(p_n^\sS (a,\mu), p_{n-1}^\sS (a, \mu )) \not\equiv
(p_n^{\sS'} (a,\mu ) , p_{n-1}^{\sS'} (a,\mu )) \,.
\eeq
\end{lemma}

\paragraph{Proof.}
It suffices to show that, given $\mu \neq 0$, 
the sequence $\sS$ is uniquely reconstructible from the data
$(p_{n} (a, \mu )$, 
$p_{n-1} (a,\mu)) := (p_n^\sS (a,\mu ), p_{n-1}^\sS (a,\mu ))$.
We define
\beql{eq88}
(p_1^\sS (a,\mu ), p_0^\sS (a,\mu )) := (1,0)
\eeq
and we have, for $1 \le j \le n$,
\beql{eq89}
\left[\begin{array}{c}
p_j^\sS (a,\mu ) \\ p_{j-1}^\sS (a,\mu )
\end{array}\right]
=
\left[ \begin{array}{cc}
F_{ab} (S_j) & -1 \\1 & 0
\end{array} \right] 
\left[\begin{array}{c}
p_{j-1}^\sS (a,\mu ) \\ p_{j-2}^\sS (a,\mu )
\end{array}\right]
\,.
\eeq
By induction $j \ge 1$, $p_j^\sS (a,\mu ) \in \RR[a]$ is a monic
 polynomial of degree $j-1$.
Applying \eqn{eq89}, given $j=n$, the given data 
$(p_n (a,\mu ) , p_{n-1} (a, \mu ))$ are 
monic polynomials of degree $n$, $n-1$ respectively.
Now there is a
unique choice of symbol $S_n$ such that
$$p_{n-2} (a, \mu ) := F_{ab} (S_n) p_{n-1} (a,\mu ) - p_n (a,\mu )$$
is a polynomial in $a$ of degree at most $n-2$.
Thus $S_n$ is uniquely determined, so we can calculate 
$(p_{n-1}^S (a,\mu ), p_{n-2}^S (a,\mu ))$.
This process can now be repeated to successively determine 
$S_{n-1}$, $S_{n-2}$, $\ldots$, $S_2, S_1$.~~~$\bsq$ \\

We next treat the case corresponding to $\mu=0$, where we obtain
univariate polynomials. 

%
%
%
\begin{lemma}\label{le81}
Let the univariate polynomials 
$P_n (a) \in \ZZ [a]$ for $n \in \ZZ$ be given by
the recurrence
\beql{eq81}
p_n (a) = a p_{n-1} (a) - p_{n-2} (a) ~,
\eeq
with the 
initial conditions $p_0 (a) =0$, $p_1 (a) =1$.
Then $p_{-n} (a) = - p_n (a)$ and, for $n \ge 2$,
\beql{eq82}
p_n (a) = \prod_{j=1}^{n-1} \left( a - \cos \frac{\pi j}{n} \right) ~.
\eeq
For $n = q \ge 3$ a prime, the irreducible factorization 
of $p_q (a)$ over $\QQ [a]$ is
\beql{eq83}
p_q (a) = (-1)^{\frac{n-1}{2}} p_q^{(0)} (a) p_q^{(0)} (-a)
\eeq
where
\beql{eq84}
p_q^{(0)} (a) =
\prod_{j=1}^{\frac{q-1}{2}} \left( a - \cos \frac{\pi j}{q} \right) 
\in \ZZ [a]
\eeq
has degree $\frac{q-1}{2}$.
\end{lemma}

\paragraph{Proof.}
The relation $p_{-n} (a) = p_n (a)$ 
is easily checked by induction on $n \ge 1$.
For $n=1,2$ we have
$$p_n \left( x+ \frac{1}{x} \right) = 
\frac{x^n - \left( \frac{1}{x} \right)^n}{x- \frac{1}{x}} ~,
$$
and this relation holds for all $n \ge 3$, by
induction on $n$, verifying the recurrence \eqn{eq81}.
Thus
\beql{eq84a}
p_n \left( x+ \frac{1}{x} \right) = x^{1-n} \frac{x^{2n}-1}{x^2-1} ~.
\eeq
The right side of \eqn{eq84a} clearly has zeros at 
$x = e^{\frac{\pi ij}{n}}$ for $1 \le j \le n-1$ and $n+1 \le j \le 2n-1$.
Now
$$a = x+ \frac{1}{x} = e^{\frac{\pi ij}{n}} +
e^{- \frac{\pi ij}{n}} = 2 \cos
\frac{\pi j}{n} ~,
$$
and these take $n-1$ distinct values for 
$1 \le j \le n-1$ (repeated for $n+1 \le j \le 2n$).
This accounts for $n-1$ distinct roots of the 
polynomial $p_n (a)$, and since it is a monic polynomial 
of degree $n-1$, the factorization \eqn{eq82} follows.
In fact $p_n (x) = U_n ( \frac{1}{2} a)$,
 where $U_n (a)$ is the Chebyshev polynomial of the second kind,
cf. Rivlin \cite[Chap. 5]{Ri90}.

Now suppose $n = q \ge 3$ is an odd prime. 
Then  we have
$\cos \frac{\pi (q-j)}{q} = - \cos \frac{\pi j}{q}$, which yields the
factorization \eqn{eq83}, \eqn{eq84}  over $\QQ [a]$.
Finally  $p^{(0)}_q (a) \in \ZZ [a]$ is irreducible 
because its roots are a complete set of Galois conjugates in 
$\QQ (\zeta_q + \zeta_q^{-1} )$, with $\zeta_q = \exp(\frac{2 \pi i}{q})$,
which is a field of degree $\frac{q-1}{2}$ over $\QQ$. $~\bsq$. \\

%
%
%
\paragraph{Proof of Theorem \ref{th82}.}
By Theorem \ref{th62}(i) a rational rotation interval 
$I_\mu (\frac{m}{n} )$ is degenerate if and only if the 
corresponding $T_{\mu\nu}$ is periodic.
By Theorem 2.3 of part I 
this occurs if and only if $(0,1)$ is a 
periodic point, with
$$
T_{\mu \nu}^{(n)}(0,1) = (0,1) ~.
$$
Since $T_{\mu\nu} (0,1) = (-1,0)$ we can write this condition as
\beql{eq810}
\left[ \begin{array}{cc}
F_{ab}(S_{n-1}) & -1 \\ 1 & 0
\end{array}
\right] \cdots \left[
\begin{array}{cc}
F_{ab}(S_1 ) & -1 \\
1 & 0
\end{array}
\right]
\left[
\begin{array}{c} -1\\0\end{array}\right]
= 
\left[
\begin{array}{c} 0\\1\end{array}\right]
~,
\eeq
where $\sS := (S_1, S_2, \ldots, S_n )$ is a certain symbol 
sequence with each $S_i = \pm 1$, and we define
$$F_{ab} (S_i) := 
\left\{
\begin{array}{lll}
a & \mbox{if} & S_i = 1 , \\
b = a-2\mu & \mbox{if} & S_i =-1.
\end{array}
\right.
$$
We define
$$(p_1^\sS (a,\mu), p_0^\sS (a,\mu)) = (1,0)$$
and set
\beql{eq811}
-
\left[
\begin{array}{c}
p_j^\sS (a,\mu)\\p_{j-1}^\sS (a,\mu)
\end{array}
\right]
=
\left[ \begin{array}{cc}
F_{ab}(S_j ) & -1 \\ 1 & 0 \end{array}\right] \cdots
\left[ \begin{array}{cc}
F_{ab}(S_1) & -1 \\ 1 & 0
\end{array}
\right]
\left[
\begin{array}{c}
-1\\0
\end{array}
\right]
 \,.
\eeq
The polynomials $p_j^\sS (a,\mu )$ satisfy the recursion
\beql{eq812}
p_j^\sS (a,\mu) := F_{ab}(S_j) p_{j-1}^\sS (a,\mu) - p_{j-1}^\sS (a,\mu ) ~,
\eeq
which was studied in Lemma \ref{le82}.
For $\mu =0$, $F_{\mu\nu}(S) \equiv a$ and the recurrence is independent 
of the symbol sequence $\sS$, and $p_j^S (a,0) = p_j (a)$, 
the polynomials in Lemma \ref{le81}.
The definition \eqn{eq811} for $j=n$ put in \eqn{eq810} becomes
$$
- (p_n^\sS (a,\mu) , p_{n-1}^\sS (a,\mu )) = (0,1) ~,
$$
so that, viewing $\mu \in \RR$ as fixed, the two polynomials

\beql{eq813}
p_n^\sS (a,\mu)=0 ~~~\mbox{and}~~~ p_{n-1}^\sS (a,\mu )+1=0
\eeq
have a common root in the $a$-variable.
For a symbol sequence $\sS$ corresponding to the 
rotation interval $I_\mu (\frac{m}{n} )$, 
a common root is $a$ is a necessary and sufficient condition for degeneracy.

We now study the common root condition 
for an arbitrary symbol sequence $\sS = (S_1, S_2, \ldots, S_{n-1} )$.
We view $\mu$ as a second indeterminate so that 
$p_n^\sS (a,\mu )$ is a bivariate polynomial in $\ZZ [a,\mu ]$.
Consider the resultant
$$
R_n^\sS (\mu ) := 
Res_a \left[p_n^\sS (a,\mu) , p_{n-1}^\sS (a,\mu) -1 \right]~.
$$
It is immediate that $R_n^\sS (\mu) \in \ZZ [\mu ]$, 
using the determinant formula for the resultant. A common
root occurs at $\mu= \mu_0$ if the resultant vanishes at $\mu_0$.
There are two cases. 

\paragraph{Case 1.}
{\em $R_n^\sS (\mu )$ is not identically  $0$.} \\

Since $R_n^\sS (\mu ) \in \ZZ [\mu ]$, all roots of 
$R_n^\sS (\mu )$ are algebraic numbers.
Also $R_n^\sS (\mu )$ is of degree at most $2n-1$, 
so there are at most $2n-1$ such roots.
For all $\mu$ not in this set, the equations \eqn{eq813} 
didn't have a common root.

\paragraph{Case 2.}
{\em $R_n^S (\mu )$ is identically $0$.} \\

In this case the two polynomials in \eqn{eq813} have a 
common bivariate factor $q(a,\mu )$, and we have
\beql{eq814}
p_n^\sS (a,\mu ) = q(a,\mu) \tilde{q} (a,\mu ) ~,
\eeq
in which $1 \le \deg_a (q(a, \mu)) \le n-1$.
By Gauss' lemma $q(a,\mu) \in \ZZ [a, \mu ]$, and is a 
monic polynomial (since $p_n^S (a,\mu )$ is monic).

We now suppose that $n=q \ge 3$ is prime.
The proof proceeds in three steps, which we will subsequently establish.

\paragraph{Step 1.}
The polynomials $q(a,\mu )$ and $\tilde{q}(a,\mu )$ are 
necessarily monic polynomials in $\ZZ[a,\mu ]$ of degree 
$\frac{p-1}{2}$ in $a$.

\paragraph{Step 2.}
A nontrivial factorization \eqn{eq814} over $\ZZ[a,\mu ]$ 
exists if and only if the sequence 
$\sS = (S_1, S_2, \ldots, S_{p-1} )$ is a palindrome, 
i.e. $S_j =S_{p-j}$ for $1 \le j \le \frac{p-1}{2}$.

\paragraph{Step 3.}
For any degenerate rotation interval $I_\mu (\frac{k}{p} )$,
$1 \le k \le p-1$, its unique legal symbol 
sequence $\sS$ has $S_1 = -1$ and $S_{p-1} =1$, 
so $\sS$ is not a palindrome. \\

Step 1 follows using Lemma \ref{le81}.
Given the bivariate
factorization \eqn{eq814} over $\ZZ[a,\mu ]$, setting $\mu=0$ gives
\beql{eq815}
p_q (a) = q(a,0) \tilde{q}(a,0) ~.
\eeq
Now $\deg_a (q(a,0)) = \deg_a (q(a,\mu ))$ because $q(a,\mu )$ 
is a monic polynomial in $a$;
similarly for $\tilde{q} (a,0)$.
Lemma \ref{le81} says that $p_q (a)$ factors over $\ZZ (a)$ into two
irreducible factors of degree $\frac{q-1}{2}$, necessarily $q(a,0)$ 
and $\tilde{q}(a,0)$ must be these two factors.
We conclude that
$$\deg_a (q(a,\mu)) = \deg_a (\tilde{q} (a,\mu )) = \frac{q-1}{2} ~,$$
completing step 1.

Step 2 will follow using Lemma \ref{le82}.
For an arbitrary $n \ge 1$ we can apply the Euclidean algorithm.
We write 
\beql{eq816}
(p_n^\sS (a,\mu ), p_{n-1}^\sS (a,\mu ) + 1) =
(p_n^\sS (a,\mu) + p_0^{\sS^R} (a,\mu ) , 
p_{n-1}^\sS (a,\mu ) + p_1^{\sS^R} (a,\mu )) ~,
\eeq
in which $\sS^R$ is the reversed sequence
$$\sS^R := (S_{n-1}, S_{n-2}, \ldots, S_1 ).$$
Now the recursion \eqn{eq812} for $\sS$ and $\sS^R$ gives
\begin{eqnarray*}
&& (p_{n-j}^\sS (a,\mu ) + p_j^{\sS^R} (a,\mu )) - 
F_{ab} (\sS_{n-j} ) (p_{n-j-1}^\sS (a,\mu) + p_{j+1}^{\sS^R} (a,\mu )) \\
&& \qquad\qquad = - (p_{n-j-2}^\sS (a,\mu) + p_{j+2}^{\sS^R} (a,\mu )) ~.
\end{eqnarray*}
We obtain, by induction on $j \ge 1$, that
\begin{eqnarray*}\label{eq817}
G(a,\mu) & = & \mbox{g.c.d.} (p_n^\sS (a,\mu) , 
p_{n-1}^\sS (a,\mu ) -1 ) \nonumber \\
& = & \mbox{g.c.d.} (p_{n-j}^\sS (a,\mu ) + p_j^{\sS^R} (a,\mu ) ,
p_{n-j-1}^\sS (a,\mu ) + p_{j+1}^{\sS^R} (a,\mu ) \,.
\end{eqnarray*}
Now suppose
$n=2m+1$ and choose $j= \frac{m-1}{2}$.
Then the right side terms are
\begin{eqnarray*}
r_1^\sS (a,\mu ) & := & p_{\frac{n+1}{2}}^\sS (a,\mu ) + 
p_{\frac{n-1}{2}}^{S^4} (a,\mu ) \\
r_2^\sS (a,\mu ) & = & p_{\frac{n-1}{2}}^\sS (a,\mu ) + 
p_{\frac{n+1}{2}}^{\sS^R} (a,\mu ) ~.
\end{eqnarray*}
These polynomials lie in $\ZZ[a,\mu ]$, 
and in the $a$-variable they are both of 
degree $\frac{n-1}{2}$, with monic top degree term $a^{\frac{n-1}{2}}$.
Thus $r^S (a,\mu ) = r_1^\sS (a,\mu ) - r_2^\sS (a,\mu )$ 
has $\deg_a (r^\sS (a,\mu )) < \frac{n-1}{2}$.

Now we impose the stronger condition that  $n=q \ge 3$ is prime.
There are two subcases.

\paragraph{Subcase 1.}
$r^\sS (a,\mu ) \not\equiv 0$. \\

Now $G(a,\mu )$ divides $r^\sS (a,\mu )$ hence we conclude
$$\deg_a (G(a,\mu )) < \frac{q-1}{2} ~.$$
But by step 1, any nontrivial divisor of $p_n^\sS (a,\mu )$ 
has degree at least $\frac{q-1}{2}$, a contradiction.
Thus Subcase 1 never occurs.

\paragraph{Subcase 2.}
$r^\sS (a,\mu ) \equiv 0$. \\

We now have
\beql{eq818}
r^\sS(a,\mu ) = (p_{\frac{q+1}{2}}^\sS (a, \mu ) - 
p_{\frac{q+1}{2}}^{\sS^R} (a,\mu ) +
( p_{\frac{q-1}{2}}^\sS (a,\mu ) - p_{\frac{q-1}{2}}^{\sS^R} (a,\mu )) = 0~.
\eeq
We assert that this equation implies the separate equalities
\begin{eqnarray}\label{eq819}
p_{\frac{q+1}{2}}^\sS (a, \mu ) & = & 
p_{\frac{q+1}{2}}^{\sS^R} (a,\mu ) ~, \nonumber \\
p_{\frac{q-1}{2}}^\sS (a,\mu ) & = & p_{\frac{q-1}{2}}^{\sS^R} (a,\mu ) ~.
\end{eqnarray}
Indeed, considering these polynomials in $\ZZ[a,\mu ]$, 
the monomials in \eqn{eq817} must cancel term-by-term.
However every monomial appearing in $p_j^\sS (a,\mu )$ 
for $j \ge 1$ has its total degree congruent to
$j-1$ $(\bmod ~2 )$, as one proves by induction on $j \ge 1$.
Thus each monomial appearing in
$p_{\frac{q+1}{2}}^\sS (a,\mu )$ can only be cancelled 
by a monomial in $p_{\frac{q+1}{2}}^{\sS^R} (a,\mu )$, 
and similarly for $p_{\frac{q-1}{2}}^\sS (a,\mu )$ and
$p_{\frac{q-1}{2}}^{\sS^R} (a,\mu )$.
Thus \eqn{eq819} follows.

We next apply Lemma \ref{le82} with $n= \frac{q+1}{2}$ 
to \eqn{eq819}, to conclude that $\sS = \sS^R$, 
so the sequence $\sS$ must be a palindrome.
Conversely, if $\sS$ is a palindrome, then \eqn{eq819} holds, whence
$r^\sS (a,\mu ) =0$, so that
$$ q(a, \mu):= \mbox{g.c.d.} (p_q^\sS (a,\mu ) , p_{q-1}^\sS (a,\mu ) -1) =
p_{\frac{q+1}{2}}^\sS (a,\mu ) - p_{\frac{q-1}{2}}^\sS (a,\mu)$$
is a nontrivial factor satisfying \eqn{eq814}.
This completes step 2.

Step 3 is a calculation. For a degenerate rotation interval we must have 
$T_{\mu v_0}^{(q)} (0,1) = (0,1)$.
The sign $S_1$ is associated to $T_{\mu v_0} (-1,0)$, hence $S_1=-1$.
The preimage $T_{\mu v_0}^{(-1)} (0,1)$ 
necessarily has the form $(1, \ast )$, hence $S_{q-1} = +1$.
Thus the associated symbol sequence $\sS$ is not a palindrome.
This completes step 3.

We conclude from steps 1--3 combined
that Case 2 can never occur for any degenerate
rotation interval
$I_\mu (\frac{j}{q} )$ where $q \ge 3$ is prime 
and $1 \le j \le \frac{q-1}{2}$.

To complete the proof of Theorem \ref{th82}, suppose that  $\mu$ 
is a value such that $I_\mu (\frac{j}{q} )$ is degenerate. 
Then $\mu$ must be  a value allowed by Case 1 for some symbol sequence 
$\sS = (S_1, \ldots, S_{q-1} )$.
There are at most $2^{q-1}$ such sequences, each of which 
has at most $2q-1$ allowed values of $\mu$, all algebraic numbers.
Thus $E(\frac{j}{q} )$ contains 
at most $2^{q-1} (2q-1)$ elements.~~~$\bsq$ \\

We are now in a position to
establish that $Spec_{\infty}[\mu]$ is a Cantor set whenever $\mu$
is a transcendental number.

%
%
%
\paragraph{Proof of Theorem \ref{th29}.}
Let $\mu \in \RR$ be fixed.
We show $Spec_{\infty}[\mu] $
 is totally disconnected  whenever
 $\mu \in \RR$ is transcendental.
The set $Spec_{\infty}[\mu] $ is closed, since it is the 
intersection of the closed sets 
$\{ (\mu, \nu ): -\infty < \nu < \infty \}$ and 
 $\Omega_{SB}$.
To show it is totally disconnected it suffices 
to show that between any two points 
$\nu_0, \nu_1 \in  Spec_{\infty}[\mu]$ there is 
an open interval in its complement.
If $\nu_0, \nu_1$ have the same rotation number they are 
necessarily the endpoints of a nontrivial rational rotation interval 
$I_{\mu_0} (r)$ whose interior is in the complement of 
$Spec_{\infty}[\mu]  $.
If the rotation numbers are unequal, say 
$r(S_{\mu_0 \nu_0} ) < r(S_{\mu_1 \nu_1})$ 
then one can find a rational $\frac{j}{q}$ with $q$ prime such that
$$r(S_{\mu_0 \nu_0} ) < \frac{j}{q} < r(S_{\mu_1 \nu_1} ) ~.$$
Theorem \ref{th82} guarantees that the rotation interval 
$I_{\mu} (\frac{j}{q} )$ is nontrivial, 
and by Theorem 2.4 of part I  its interior is not in $\Omega_{SB}$.
Thus $Spec_{\infty}[\mu] $ is a totally disconnected set.

The set $ Spec_{\infty}[\mu]  $ is a perfect set for all $\mu \in \RR$.
Indeed, given $\nu \in Spec_{\infty}[\mu]$, 
the rotation number $r(S_{\mu \nu} )$ can be approximated 
both from above and below by sequences
$(\mu,\nu_n^{+}),  (\mu, \nu_n^{-}) \in \Omega_{OB}$ having 
irrational rotation numbers
converging to $r(S_{\mu \nu} )$.
Then at least one of 
$\lim_{n \to \infty} \nu_n^- = \nu$ or 
$\lim_{n \to \infty} \nu_n^+ = \nu$ holds, 
and both hold if $r(S_{\mu \nu} )$ is irrational.
Thus $ Spec_{\infty}[\mu]$ is a perfect set.~~~$\bsq$

%
%
%
%
%

\section{Concluding Remarks}

Does there exist a value $\mu \neq 0$ 
for which $Spec_{\infty}[\mu]$ is {\em not} a Cantor set?
For algebraic $\mu \neq 0$ extra degeneracies 
for rational rotation intervals can occur, but we do not know 
whether there are ever enough of them to destroy 
the property of $Spec_{\infty}[\mu]$ being totally disconnected.
We numerically considered the test case $\mu_0 = \frac{1}{2}$.
Its full list of degenerate rotation intervals 
with denominator below $64$ is given by 
$r = \frac{1}{5}$, $\frac{2}{7}$, $\frac{2}{15}$, 
$\frac{7}{20}$ and $\frac{11}{27}$, plus members of the sequence 
$r= \{ \frac{2n-1}{4n} : n \ge 2 \}$ given in Example 3.3 of part I.
This list shows no sign of covering all rationals sufficiently near
one particular point, as would be required if 
$Spec_{\infty}[\mu]  $ were not to be  totally disconnected.

     Another open  question concerns whether the set of $\nu$
such that $T_{\mu\nu}$ has a nontrivial bounded orbit is of positive 
one-dimensional Lebesgue measure.
The only case for which the answer is known is $\mu=0$, which gives measure
$4$ (the maximum possible for any $\mu$).  This is a very
special case, since for $\mu=0$, all rational rotation number intervals
have length 0.  For $\mu \neq 0$,
Theorem~\ref{th82} shows that  infinitely many of the rational
rotation number intervals have positive length, so the uncovered
measure is strictly smaller than $4$.
Numerical experiments 
suggest, in two different ways,
that the measure is positive for at least some non-zero $\mu$.
First, the calculation at $\mu=1$ of the uncovered area
as a function of denominator suggested positive measure
will remain. Second, numerical random selection of 
parameter values in the allowed interval seemed to produce invariant circles
with positive probability. \\


\end{document}